\documentclass[12pt]{amsart}
\usepackage{amssymb}
\usepackage{graphics}
\usepackage{latexsym}
\usepackage{amsmath}
\usepackage{amssymb,amsthm,amsfonts}
\usepackage{amscd}
\usepackage[arrow, matrix, curve]{xy}
\usepackage{syntonly}
\title[Arakelov Inequality and $p$-Rank Zero Locus]{An Arakelov Inequality in Characteristic $p$ and Upper Bound of $p$-Rank Zero Locus}
\author[Jun Lu]{Jun Lu}
\author[Mao Sheng]{Mao Sheng}
\author[Kang Zuo]{Kang Zuo}
\address{Universit\"{a}t Mainz, Fachbereich 17, Mathematik, 55099 Mainz, Germany}
\email{lu@mathematik.uni-mainz.de}
\address{Department of Mathematics, East China Normal
University, 200062 Shanghhai, P.R. China}
\email{msheng@math.ecnu.edu.cn}
\address{Universit\"{a}t Mainz, Fachbereich 17, Mathematik, 55099 Mainz, Germany}
\email{kzuo@mathematik.uni-mainz.de}

\begin{document}
\theoremstyle{plain}
\newtheorem{thm}{Theorem}[section]
\newtheorem{theorem}[thm]{Theorem}
\newtheorem{lemma}[thm]{Lemma}
\newtheorem{corollary}[thm]{Corollary}
\newtheorem{proposition}[thm]{Proposition}
\newtheorem{addendum}[thm]{Addendum}
\newtheorem{variant}[thm]{Variant}
\theoremstyle{definition}
\newtheorem{construction}[thm]{Construction}
\newtheorem{notations}[thm]{Notations}
\newtheorem{question}[thm]{Question}
\newtheorem{problem}[thm]{Problem}
\newtheorem{remark}[thm]{Remark}
\newtheorem{remarks}[thm]{Remarks}
\newtheorem{definition}[thm]{Definition}
\newtheorem{claim}[thm]{Claim}
\newtheorem{assumption}[thm]{Assumption}
\newtheorem{assumptions}[thm]{Assumptions}
\newtheorem{properties}[thm]{Properties}
\newtheorem{example}[thm]{Example}
\newtheorem{conjecture}[thm]{Conjecture}
\numberwithin{equation}{thm}

\newcommand{\pP}{{\mathfrak p}}
\newcommand{\sA}{{\mathcal A}}
\newcommand{\sB}{{\mathcal B}}
\newcommand{\sC}{{\mathcal C}}
\newcommand{\sD}{{\mathcal D}}
\newcommand{\sE}{{\mathcal E}}
\newcommand{\sF}{{\mathcal F}}
\newcommand{\sG}{{\mathcal G}}
\newcommand{\sH}{{\mathcal H}}
\newcommand{\sI}{{\mathcal I}}
\newcommand{\sJ}{{\mathcal J}}
\newcommand{\sK}{{\mathcal K}}
\newcommand{\sL}{{\mathcal L}}
\newcommand{\sM}{{\mathcal M}}
\newcommand{\sN}{{\mathcal N}}
\newcommand{\sO}{{\mathcal O}}
\newcommand{\sP}{{\mathcal P}}
\newcommand{\sQ}{{\mathcal Q}}
\newcommand{\sR}{{\mathcal R}}
\newcommand{\sS}{{\mathcal S}}
\newcommand{\sT}{{\mathcal T}}
\newcommand{\sU}{{\mathcal U}}
\newcommand{\sV}{{\mathcal V}}
\newcommand{\sW}{{\mathcal W}}
\newcommand{\sX}{{\mathcal X}}
\newcommand{\sY}{{\mathcal Y}}
\newcommand{\sZ}{{\mathcal Z}}
\newcommand{\A}{{\mathbb A}}
\newcommand{\B}{{\mathbb B}}
\newcommand{\C}{{\mathbb C}}
\newcommand{\D}{{\mathbb D}}
\newcommand{\E}{{\mathbb E}}
\newcommand{\F}{{\mathbb F}}
\newcommand{\G}{{\mathbb G}}
\newcommand{\HH}{{\mathbb H}}
\newcommand{\I}{{\mathbb I}}
\newcommand{\J}{{\mathbb J}}
\renewcommand{\L}{{\mathbb L}}
\newcommand{\M}{{\mathbb M}}
\newcommand{\N}{{\mathbb N}}
\renewcommand{\P}{{\mathbb P}}
\newcommand{\Q}{{\mathbb Q}}
\newcommand{\R}{{\mathbb R}}
\newcommand{\SSS}{{\mathbb S}}
\newcommand{\T}{{\mathbb T}}
\newcommand{\U}{{\mathbb U}}
\newcommand{\V}{{\mathbb V}}
\newcommand{\W}{{\mathbb W}}
\newcommand{\X}{{\mathbb X}}
\newcommand{\Y}{{\mathbb Y}}
\newcommand{\Z}{{\mathbb Z}}
\newcommand{\id}{{\rm id}}
\newcommand{\rank}{{\rm rank}}
\newcommand{\END}{{\mathbb E}{\rm nd}}
\newcommand{\End}{{\rm End}}
\newcommand{\Hom}{{\rm Hom}}
\newcommand{\Hg}{{\rm Hg}}
\newcommand{\tr}{{\rm tr}}
\newcommand{\Sl}{{\rm Sl}}
\newcommand{\Gl}{{\rm Gl}}
\newcommand{\Cor}{{\rm Cor}}
\newcommand{\proofend}{\hspace*{13cm} $\square$ \\}
\maketitle
\begin{abstract}
In this paper we show an Arakelov inequality for semi-stable
families of algebraic curves of genus $g\geq 1$ over characteristic
$p$ with nontrivial Kodaira-Spencer maps. We apply this inequality
to obtain an upper bound of the number of algebraic curves of
$p-$rank zero in a semi-stable family over characteristic $p$ with
nontrivial Kodaira-Spencer map in terms of the genus of a general
closed fiber, the genus of the base curve and  the number of
singular fibres. An extension of the above results to smooth
families of Abelian varieties over $k$ with $W_2$-lifting assumption
is also included.
\end{abstract}

\footnotetext[1]{This work was supported by the SFB/TR 45 ¡°Periods,
Moduli Spaces and Arithmetic of Algebraic Varieties¡± of the DFG
(German Research Foundation).}

\footnotetext[2]{The second named author is supported by a
Postdoctoral Fellowship in the East China Normal University.}

\section{Introduction}
Let $f: X\to C$ be a non-isotrivial semi-stable family of algebraic
curves of genus $g \geq 1$ over smooth projective curve $C$ over
complex numbers. Let $S\subset C$ be the singular locus over which
the fibration $f$ degenerates. The classical Arakelov inequality
(cf. \cite{Fa}) states that the following inequality holds
$$
\deg f_{*}\omega_{X/C}\leq \frac{g}{2}\deg \Omega_{C}(S).
$$
This is one of key ingredients in the proof by Arakelov on the
Shararevich conjecture that the isomorphism classes of genus $g$
curves over a given functional field with fixed degeneracy are
finite. There are new developments since then on the Arakelov
inequalities and their applications to certain geometric problems on
moduli spaces of polarized algebraic manifolds. (For more
information, refer to the recent survey articles \cite{V}, \cite{Z}
and references therein.) In a series of papers \cite{VZ01},
\cite{VZ06} and \cite{MVZ07}, the following generalized form of
Arakelov inequality was obtained:
\begin{theorem}[Theorem 1.1, \cite{V}; Theorem 4.4, \cite{Z}, $n=1$
case]\label{Arakelov inequality for families of curves in
characteristic 0} Let $f: X\to C$ be a semi-stable family of curves
as above. Let $\sE$ be a coherent subsheaf of
$f_{*}\omega_{X/C}^{\nu},\ \nu\geq 1$. Then the following inequality
holds:
$$
\frac{\deg\sE}{\rank \sE}\leq \frac{\nu}{2} \deg \Omega_{C}(S).
$$
\end{theorem}
In this paper we give a characteristic $p$ analogue of the Arakelov
inequality in the above generalized form. Let $k$ be the algebraic
closure of finite field $\F_{p}$ with $p$ an odd prime. Now we let
$f: X\to C$ be a semi-stable family of algebraic curves of genus
$g\geq 1$ with nontrivial Kodaira-Spencer map over projective curve
$C$, which is defined over $k$. The nontriviality of Kodaira-Spencer
map means that it is nonzero at one closed point in the smooth locus
of the base curve and this assumption is equivalent to saying that
the family is non-isotrivial and it is not semi-stable reduction of
the base change of another family $f': X'\to C$ under the Frobenius
map $F_{C}: C\to C$.
\begin{theorem}[Theorem \ref{Arakelov inequality for families of curves in characteristic
p}] Let $f: X\to C$ be a semi-stable family of algebraic curves of
genus $g\geq 1$ over $k$ with nontrivial Kodaira-Spencer map. Let
$\sE$ be a coherent subsheaf of $f_{*}\omega_{X/C}^{\nu},\ \nu\geq
1$. Then the following strict inequality holds:
\begin{eqnarray*}
\frac{\deg \sE}{\rank \sE}< 2\nu g\deg \Omega_{C}(S).
\end{eqnarray*}
\end{theorem}
One can not deduce the above Arakelov inequality in characteristic
$p$ directly from the Arakelov inequality in characteristic zero as
given in Theorem \ref{Arakelov inequality for families of curves in
characteristic 0}. This is because there exists non-liftable family
of algebraic curves in characteristic $p$, and even if the family
$f$ in characteristic $p$ comes from the reduction of family at a
prime over $p$, there exists possibly non-liftable coherent subsheaf
in $f_{*}\omega_{X/C}^{\nu}$. Moreover, in characteristic 0 G.
Faltings \cite{Fa} showed the Arakelov inequality for a semi-stable
family of Abelian varieties and one can deduce the classical
Arakelov inequality for family of algebraic curves by considering
the associated family of Jacobians. In characteristic $p$ we do not
have the same form of Arakelov inequality for semi-stable families
of Abelian varieties. Actually Moret-Bailly \cite{Mo} constructed a
semi-stable family of genus 2 curves over $\P^1$ defined over $k$,
whose associated Jacobian fibration is a smooth family of
supersingular Abelian surfaces. It was shown furthermore that
$f_{*}\omega_{X/C}=\sO_{\P^1}(p)\oplus \sO_{\P^1}(-1)$ (See 3.2,
\cite{Mo}). It seems that a properly formulated Arakelov inequality
for semi-stable families of Abelian varieties in characteristic $p$
has not been known yet.\\

Besides giving a characteristic $p$ analogue of the classical
Arakelov inequality, another motivation for us is using it to study
certain geometric problems of the moduli spaces of curves defined
over $k$. Please read the articles by F. Oort \cite{Oort}, G. van
der Geer, or both of them in \cite{FL} for a general introduction to
the moduli spaces of curves and Abelian varieties in characteristic
$p$. Particulary we shall use the same notations as them and shall
not repeat the definitions if they have already been defined there.
We start our discussions with a very classical example, the Legendre
family of elliptic curves over $k$. It is given by the affine
equation
$$
y^2=x(x-1)(x-t)
$$
with parameter $t$. As well known, the family has three singular
fibers and $\frac{p-1}{2}$ (smooth) supersingular fibers. One needs
to do semi-stable reduction to the Legendre family in order to fit
into our considerations. A possible base change is given by the
double cover of the base curve $\P^1$ branching at $0$ and $\infty$.
The obtained semi-stable model of the Legendre family has then four
singular fibers and $p-1$ supersingular fibers. The base curve is
actually reduction of the modular curve with level structure
$\Gamma(2)\cap \Gamma_{1}(4)\subset SL(2,\Z)$. This is one of six
examples of semi-stable families of elliptic curves over $\P^1$ due
to A. Beauville \cite{Ar}. From these information plus that of
automorphisms of elliptic curves one can deduce the classical
Deuring mass formula on the number of supersingular elliptic curves
in the coarse moduli space over $k$. It is remarkable that T.
Ekedahl and G. van der Geer have generalized the formula to the
Ekedahl-Oort strata in the coarse moduli space of principally
polarized Abelian varieties over $k$ (see the article by G. van der
Geer in \cite{FL} and references therein). One can ask as next step
for formulas or equalities for Ekedahl-Oort strata in the Torelli
locus. Unfortunately, such questions remain difficult in general
(see the article by F. Oort in \cite{FL}). In this paper we want to
use the above Arakelov inequality to present certain inequality
about $p$-rank zero locus. To motivate it, it can be shown that for
example $p-1$ is the maximal number of supersingular elliptic curves
in a semi-stable family of elliptic curves over $\P^1$ with four
singular fibers over $k$, whose Kodaira-Spencer map is nontrivial.
(See discussions after Proposition \ref{elliptic curve case}.) For a
higher genus semi-stable fibration we show the following strict
inequality:
\begin{theorem}[Theorem \ref{upper bound of p-rank zero locus}]
Let $f: X\to C$ be a semi-stable family of algebraic curves of genus
$g\geq 1$ over smooth projective algebraic curves $C$ over $k$ with
nontrivial Kodaira-Spencer map. If the $p$-rank of the generic fiber
of $f$ is nonzero, then the number of $p$-rank zero closed fibers of
$f$ is bounded from above strictly by
$$[2p^{g}(2g^2-1)+2g(g-1)](2b-2+s)$$
where $b$ is the genus of base curve $C$ and $s$ is the number of
singular fibers of $f$.
\end{theorem}
The contents of the paper are organized as follows. In section 2 we
prove the Arakelov inequality in characteristic $p$. In section 3 we
discuss the relative Frobenius morphism of $f: X\to C$ and show an
inequality for the slopes of coherent subsheaves in
$(F_{C}^*)^{n}R^{1}f_{*}\sO_{X}, \ n\geq 1$. This section serves as
preparation for section 4, but is separated from section 4 for its
own interests. In section 4 we prove the claimed upper bound of
$p$-rank zero locus.\\

\textbf{Acknowledgements:} The authors would like to thank Professor
Xiaotao Sun for useful discussions on this paper, in particular the
explanation of his results in \cite{Sun}. We would like also to
thank Professor Eckart Viehweg for discussions on the Bertnini's
theorem in characteristic $p$ and the Szpiro's inequality. Last but
not least we thank Dr. Ralf Gerkmann and Dr. Jiajin Zhang for many
general discussions and particularly on the notion of $F$-crystal.\\

\textbf{Notations and Conventions:} In the following sections the
notation $f: X\to C$ means always a semi-stable family of algebraic
curves of genus $g\geq 1$ over $k$ with nontrivial Kodaira-Spencer
map, where the base curve $C$ is smooth and projective. The set
$S\subset C$ is the singular locus of $f$. We denote by
$E^{1,0}=f_{*}\omega_{X/C}$ and $E^{0,1}=R^{1}f_{*}\sO_{X}$. They
are called to be the first and respectively the second Hodge bundles
of the family $f$, which are dual to each other by the relative
Serre duality. The slope of a coherent sheaf $\sF$ over $C$ is
defined to be $\mu(\sF)=\frac{\deg\sF}{\rank \sF}$. The $p$-rank of
a smooth projective algebraic curve of genus $g\geq 1$ is defined to
be the $p$-rank of its Jacobian (see the article by F. Oort in
\cite{FL}, or section 4 \cite{Mu} for the notion of $p$-rank). For
an algebraic variety $X$ defined over $k$, the map $F_{X}: X\to X$
is denoted to be the absolute Frobenius morphism defined by power
$p$ map on the structure sheaf of rings $\sO_{X}$ (see \cite{EV}
section 9). In section 3 we shall use occasionally the notion of
$F$-crystal and crystalline cohomology. The basic reference for
$F$-crystal is the article by N. Katz \cite{Ka} and for crystalline
cohomology is the book \cite{BO}. The rest notions on algebraic
varieties we used in the article are all standard and one can find
them in \cite{Ha}.

\section{Arakelov Inequality of Semi-stable Families of Algebraic Curves in Characteristic $p$}
Let $f: X\to C$ be a semi-stable family of algebraic curves of genus
$g\geq 1$ with nontrivial Kodaira-Spencer map. In this section we
shall prove an Arakelov inequality for the family $f$. Our proof is
based on certain techniques of algebraic surface, for which one can
consult the book by G. Xiao \cite{Xi}, and the results of L. Szpiro
in \cite{Sz}. We are going to prove the following
\begin{theorem}\label{Arakelov inequality for families of curves in characteristic p}
Let $f: X\to C$ be a semi-stable family of algebraic curves of genus
$\geq 1$ over $k$ with nontrivial Kodaira-Spencer map. For a
coherent subsheaf $\mathcal{E}$ in $f_{*}\omega_{X/C}^{ \nu},\
\nu\geq 1$ one has the following upper bound on the slope of $\sE$:
\begin{eqnarray*}\label{eq02}
\mu(\mathcal{E})< 2\nu g\deg \Omega_{C}(S).
\end{eqnarray*}
\end{theorem}

{\itshape Proof:} We proceed it by considering $g=1$ and $g\geq 2$
separately. First we study $g=1$ case. By cup product with first
Hodge bundle $E^{1,0}$ and then contraction on coefficients, the
Kodaira-Spencer map of $f$ induces a non-trivial morphism
$$
\theta: E^{1,0}\to E^{0,1}\otimes \Omega_{C}(S).
$$
Since they are invertible sheaves, $\theta$ induces an embedding
$$
(E^{1,0})^{\otimes 2}\to  \Omega_{C}(S).
$$
This implies that
$$
\deg E^{1,0} \leq \frac{1}{2}\deg \Omega_{C}(S).
$$
Because $f$ is a semi-stable elliptic fibration, one has isomorphism
$$
f_{*}\omega_{X/C}^{\nu}\simeq (f_{*}\omega_{X/C})^{\nu}.
$$
Then for any coherent subsheaf $\sE$ of $f_{*}\omega_{X/C}^{\nu}$,
which is invertible in this case, one has the following inequality
$$
\mu(\sE)=\deg \sE \leq \nu \deg E^{1,0}\leq \frac{\nu}{2}\deg
\Omega_{C}(S)< 2\nu\deg \Omega_{C}(S).
$$
Now we consider $g\geq 2$ case. Let
\begin{align*}
f_{*}\omega_{X/C}^{\nu}
=\mathcal{E}_m\supseteq\mathcal{E}_{m-1}\supseteq\cdots\supseteq\mathcal{E}_{1}\supseteq\mathcal{E}_{0}=0
\end{align*}
be the Harder-Narasimhan filtration of $f_{*}\omega_{X/C}^{\nu}$.
Obviously, it suffices to show our inequality for $\sE_1$, which has
maximal slope among all coherent subsheaves in
$f_{*}\omega_{X/C}^{\nu}$. One notes that the image sheaf of the
natural morphism $\alpha: f^{*}\sE_1\to \omega_{X/C}^{\otimes\nu}$
can be expressed by $\mathcal{I}_{Z}(\nu K_{X/C}-D)$, where $\sI_Z$
is the ideal sheaf of a zero-dimensional subscheme $Z$, $K_{X/C}$ a
relative canonical divisor and $D$ is an effective divisor. Now let
$H$ be an ample divisor of $X$.  Since $\omega_{X/C}$ is big and nef
by L. Szpiro (see Theorem 1 and Proposition 3 in \cite{Sz}), $m\nu
K_{X/C}+H$ is ample for $m\geq 0$ by Nakai's criterion (See Theorem
1 in \cite{Kl}). By Bertini's theorem, one can find a smooth curve
$$
\Gamma\in |n(m\nu K_{X/C}+H)|
$$
for $n$ large enough such that $\Gamma$ intersects transversally
with a fixed smooth closed fiber $F$ of $f$, and $\Gamma$ is not
contained in the support the kernel of $\alpha$. We put
$\pi=f|_{\Gamma}$. By construction, $\pi$ is a separable finite
morphism and the restriction of $\alpha$ to $\Gamma$
$$
\alpha|_{\Gamma}: \pi^{*}\sE_1 \to \mathcal{O}_{\Gamma}(\nu
K_{X/C}-D)
$$
is non-trivial. Because $\sE_1$ is semi-stable and $\pi$ is
separable, $\pi^{*}\mathcal{E}_1$ is still semi-stable by Lemma 3.1
\cite{Sun}. Hence we have
\begin{align*}
\mu(\pi^{*}(\mathcal{E}_1))\leq \deg\mathcal{O}_{\Gamma}(\nu
K_{X/C}-D)=(\nu K_{X/C}-D)\Gamma.
\end{align*}
We put $N=\nu K_{X/C}-D-\mu(\mathcal{E}_1)F$. Then
\begin{eqnarray*}
   N\Gamma&=&  (\nu K_{X/C}-D)\Gamma-\mu(\sE_1)F\Gamma \\
    &=& (\nu K_{X/C}-D)\Gamma-\mu(\pi^{*}(\mathcal{E}_1)) \\
    &\geq&0.
\end{eqnarray*}
It implies that $N(m\nu K_{X/C}+H)\geq 0$. Since $m$ can be
arbitrarily large, we must have $NK_{X/C}\geq 0$. Therefore, we see
that
\begin{eqnarray*}\label{mumax01}
  \nu K_{X/C}^2 &=&  (N+D+\mu(\mathcal{E}_1)F)K_{X/C} \\
   &\geq & \mu(\mathcal{E}_1)FK_{X/C} \\
    &=& (2g-2)\mu(\mathcal{E}_1).
\end{eqnarray*}
Now we recall the Szpiro's inequality (see Proposition 4.2 in
\cite{Sz})
\begin{align*}
K_{X/C}^2<(4g-4)g\deg \Omega_{C}(S).
\end{align*}
By combining the last two inequalities, we obtain therefore
\begin{align*}
\mu(\mathcal{E}_1)< 2\nu g\deg \Omega_{C}(S).
\end{align*}

\proofend

\section{Relative Frobenius Morphism and Upper Bound of Slopes of Coherent Subsheaves in the Pull-back of Second Hodge Bundle Under Iterated Frobenius Morphism}
For the family $f: X\to C$ over $k$ one has the following
commutative diagram induced by the absolute Frobenius morphisms of
the base $C$ and the total space $X$:
$$
\xymatrix{ X\ar[dr]_{f}\ar[r]^{F_{rel}}&   X^{'} \ar[d]^{f^{'}}
\ar[r]^{Pr}
                & X \ar[d]^{f}  \\
&  C \ar[r]_{F_{C}}
                & C             }
$$
where $f^{'}:X^{'}\to C$ is the base change of $f$ under $F_C$ and
$F_{rel}: X\to X^{'}$ is the so-called relative Frobenius morphism.
The $C$-morphism $F_{rel}$ is the main interest of the first half of
this section. The following simple lemma is well known. For the
convenience of the reader we include it here with a proof.
\begin{lemma}\label{Frobenius on the second Hodge bundle}
The relative Frobenius morphism $F_{rel}$ induces a natural morphism
of $\sO_{C}$-modules
$$
F_{rel}^{*}: F_{C}^{*}E^{0,1}\to E^{0,1}.
$$
\end{lemma}
{\itshape Proof:} The relative Frobenius morphism $F_{rel}$ gives
the $\sO_{X^{'}}$-morphism of sheaves
$$
\sO_{X^{'}}\to F_{rel*}\sO_{X}.
$$
It induces the morphism on the direct images
$$
R^{1}f^{'}_{*}\sO_{X^{'}}\to R^{1}f^{'}_{*}(F_{rel*}\sO_{X}).
$$
The spectral sequence of composed maps gives the natural morphism
$$
R^{1}f^{'}_{*}(F_{rel*}\sO_{X}) \to R^{1}(f^{'}\circ
F_{rel})_{*}\sO_{X}.
$$
Composition of the above two morphisms yields the morphism
$$
R^{1}f^{'}_{*}\sO_{X^{'}}\to R^{1}f_{*}\sO_{X}.
$$
Finally, by the flat base change theorem(cf. Proposition 9.3
\cite{Ha}) the $R^{1}f^{'}_{*}\sO_{X^{'}}$ is isomorphic to
$F_C^{*}R^{1}f_{*}\sO_{X}$. The lemma follows.

\proofend

In the following text we denote for brevity the $n$-th iterated
Frobenius morphism $(F_{C}^{*})^{n}$ by $F_{n}^*$. Following Lemma
\ref{Frobenius on the second Hodge bundle} one can then consider a
sequence of morphisms
$$
\Phi_{n}: F_{n}^*E^{0,1}\to E^{0,1},\  n\geq 1,
$$ which is the
composition of the following morphisms:
$$
F_{n}^*E^{0,1}\stackrel{F_{n-1}^*F_{rel}^{*}}{\longrightarrow}
F_{n-1}^*E^{0,1}\stackrel{F_{n-2}^*F_{rel}^{*}}{\longrightarrow}
\cdots \to F_{1}^{*}E^{0,1}\stackrel{F_{rel}^{*}}{\longrightarrow}
E^{0,1}.
$$
The following result describes some relations between the morphism
$\Phi_{n}$ and the $p$-rank of fibers.
\begin{proposition}\label{criteria}
Let $f: X\to C$ be a semi-stable family of algebraic curves. Let $t$
be a closed point of $C$ and the closed fiber $X_{t}$ of $f$ over
$t$ be smooth. Then the following statements hold:
\begin{itemize}
    \item [(i).] The morphism $\Phi_n$ is isomorphism at $t$ for one $n$ (or for all $n$) if and only if $X_{t}$ is
ordinary.
    \item [(ii).] The morphism $\Phi_g$ is zero at $t$ if and only if $X_{t}$ is of $p$-rank
zero.
    \item [(iii).] The morphism $\Phi_1$ is zero at $t$ only if $X_{t}$ is
supersingular.
\end{itemize}
\end{proposition}
These statements should be more or less well known. For completeness
we present a proof here, although we are aware of the possibility
that it is worse than the original one. We remark that the if-part
of (iii) above is not true.\\

They are point-wise statements. So we consider the Jacobian of an
algebraic curve of genus $g\geq 1$ or generally a principally
polarized Abelian variety $A$ of dimension $g$ over $k$. Let
$H^1_{crys}(A/W(k))$ be the first crystalline cohomology. It is
torsion free $W(k)$-module of rank $2g$. More importantly, it has a
semi-linear endomorphism $F$ on $H^1_{crys}(A/W(k))$ which makes the
pair $(H^1_{crys}(A/W(k)),F)$ the structure of $F$-crystal. By
comparison theorem, the reduction at $k$ of the $F$-crystal is
isomorphic to the Frobenius morphism on the de-Rham cohomology
$H^{1}_{dR}(A/k)$. By using the Dieudonne theory, the $p$-rank of
$A$ is equal to the multiplicity of zero slopes of the $F$-crystal.
The basic linear algebra of $F$-crystal is established in \cite{Ka}.
We shall also use a result due to N. Nygaard \cite{Ny}.\\

{\itshape Proof of Proposition \ref{criteria}:} Let $t$ be a
$k$-point of $C$, and $A=Jac(X_{t})$ be the Jacobian of $X_t$. It is
equivalent to prove the corresponding statements for
$$
(F_{A}^{*})^n: H^{1}(A,\sO_A)\to H^{1}(A,\sO_A).
$$
We look at the $F$-crystal of $A$. By the Hodge-Newton decomposition
(cf. Theorem 1.6.1 \cite{Ka}), there is a unique decomposition of
$F$-crystals
$$
(H^1_{crys}(A/W(k)),F)=(M_{0}\oplus M_{>0},F)
$$
where $(M_{0},F|_{M_{0}}$ is the unit subcrystal of
$(H^1_{crys}(A/W(k)),F)$ and $M_{>0}$ is the complement whose Newton
slopes are all positive. So the proofs of (i) and (ii) are reduced
to the following
\begin{claim}\label{relation between Frobenius on Cristalline and Frobenius on Hodge
Bundle} Let $A$ be an $g$-dimensional Abelian variety defined over
$k$. The Newton slopes of $A$ are all positive if and only if the
morphism
$$
(F_{A}^{*})^{g}: H^{1}(A,\sO_A)\to H^{1}(A,\sO_A)
$$
is zero map.
\end{claim}
{\itshape Proof of Claim:} It is known that the $k$-vector space
$H^{1}(A,\sO_A)$ is a natural quotient of $H^1_{cris}(A,W(k))$
modulo $p$, and the map $F$ modulo $p$ induces a natural morphism on
$H^{1}(A,\sO_A)$, which is identical to $F_{A}^{*}$. By N. Katz
1.3.3 \cite{Ka}, all Newton slopes are positive if and only if
$F^{2g}(H^1_{cris}(A,W(k)))\subset pH^1_{cris}(A,W(k))$, which
implies in particular
$$
(F_{A}^{*})^{2g}: H^{1}(A,\sO_A)\to H^{1}(A,\sO_A)
$$
is zero map. But since $\dim_{k}H^{1}(A,\sO_A)=g$, the $g$ times
iterate $(F_{A}^{*})^{g}$ is already zero. Conversely, if
$(F_{A}^{*})^{g}$ is a zero morphism, then
$$
(F_{A}^{*})^{g+1}: H^{1}_{dR}(A/k)\to H^{1}_{dR}(A/k)
$$
is zero. Since $2g\geq g+1$, it implies that
$$
F^{2g}(H^1_{cris}(A,W(k)))\subset pH^1_{cris}(A,W(k)).
$$
Then still by N. Katz's remark, the Newton slopes of
$H^1_{cris}(A,W(k))$ are all positive. The claims is proved.

\proofend

To the part (iii) we invoke a result of N. Nygaard.
\begin{theorem}[Nygaard, Theorem 1.2 \cite{Ny}]\label{Nygaard's
result} Let $A$ and $(H^1_{crys}(A/W(k)),F)$ be as above. Then $A$
is supersingular if and only if
\begin{itemize}
  \item [a)] $F^{g^2-g+2}$ is divisible by $p^{\frac{g^2+1}{2}-(g-1)}$ if $g$ is
  odd;
  \item [b)] $F^{g^2-2g+3}$ is divisible by
$p^{\frac{g^2+1}{2}-\frac{3}{2}(g-1)}$ if $g$ is even.
\end{itemize}
\end{theorem}
So the assumption of (iii) says that the map $F_{A}^{*}:
H^{1}(A,\sO_A)\to H^{1}(A,\sO_A)$ is zero. It implies that $F^2:
H^1_{crys}(A/W(k)) \to H^1_{crys}(A/W(k))$ is divisible by $p$. An
elementary calculation shows that the divisibility condition in
above theorem is satisfied and (iii) follows therefore from
Theorem \ref{Nygaard's result}.\\

\proofend

In the remainder of this section we shall prove the following
theorem about various upper bounds of the slopes of coherent
subsheaves contained in $F_{n}^*E^{0,1},\ n\geq 1$.
\begin{theorem}\label{Arakelov inequality for second Hodge bundle}
Let $f:X\to C$ be a semi-stable family of algebraic curves as above.
Let $\sF$ be a coherent subsheaf of $F_{n}^*E^{0,1}$ for certain
$n\geq 1$. Then we have the following upper bound of the slope of
$\sF$:
\begin{itemize}
    \item [(i).] If $f$ has an ordinary closed fiber, then $\mu(\sF)\leq
    0$;
    \item [(ii).] If $n=1$, then $\mu(\sF)< 2g(g-1)\deg
    \Omega_{C}(S)$;
    \item [(iii).] In any case one has
    $$
    \mu(\sF)< 2g(g-1)(2b-2+s)p^{n}+2(g-1)(b-1)\frac{p^n-1}{p-1},
    $$
where $b$ is the genus of $C$ and $s$ the number of singular fibers
of $f$.
\end{itemize}
\end{theorem}

{\itshape Proof of (i) and (ii):} We prove the first statement by
contradiction. By assumption that $f$ has an ordinary closed fiber,
the generic fiber of $f$ is then ordinary. By Lemma \ref{criteria}
(i), the morphism $\Phi_1$ is generically isomorphism. Now one
assumes that the coherent subsheaf $\sF$ of $F_{n}^*E^{0,1}$ has
positive slope. One can further assume that $\sF$ is of maximal
slope in $F_{n}^*E^{0,1}$. We consider the morphism
$$
F_{n}^*\Phi_{1}: F_{n+1}^*E^{0,1}\to F_{n}^*E^{0,1}
$$
and claim that the subsheaf $F_{n}^*\Phi_{1}(F_{1}^*\sF)$ of
$F_{n}^*E^{0,1}$ has larger slope than $\sF$. First of all, the rank
of $F_{n}^*\Phi_{1}(F_{1}^*\sF)$ is equal to that of $\sF$ because
$\Phi_1$ is generically isomorphism and so is $F_{n}^*\Phi_{1}$.
Next one has $$\deg F_{n}^*\Phi_{1}(F_{1}^*\sF)\geq \deg
F_{1}^*\sF=p\deg \sF.$$ Since $\deg \sF>0$, it is clear $\deg
F_{n}^*\Phi_{1}(F_{1}^*\sF)>\deg \sF$. A contradiction. \\

Now we let $\sF$ be a rank $r$ coherent subsheaf of $E^{0,1}$. The
dual of the quotient sheaf $\displaystyle E^{0,1}/\sF$ is a coherent
subsheaf of $E^{1,0}$. Furthermore, because the family $f$ has
nontrivial Kodaira-Spencer map, one has $\deg E^{1,0}>0$ by
Proposition 3, \cite{Sz}. Applying the Arakelov inequality
\ref{Arakelov inequality for families of curves in characteristic p}
in $\nu=1$ case, we have
\begin{eqnarray*}
\deg(\sF)    &=& -\deg{E^{1,0}}+\deg{(E^{0,1}/\sF)^{*}}   \\
   &< & \deg{(E^{0,1}/\sF)^{*}}\\
    &<& 2g(g-r)\deg \Omega_{C}(S).
\end{eqnarray*}
Therefore
$$
\mu(\sF)< 2g(\frac{g-r}{r})\deg \Omega_{C}(S)\leq 2g(g-1)\deg
\Omega_{C}(S).
$$
The second statement is proved.\\

\proofend

Before giving the upper bound in the general case, we shall discuss
a peculiar phenomenon about the instability of vector bundles under
Frobenius pull-back. For a coherent sheaf $\sE$ the notation
$\mu_{max}(\sE)$ (resp. $\mu_{min}(\sE)$) means the maximal (resp.
minimal) slope of coherent subsheaves in $\sE$. In \cite{Sun} X. Sun
has shown the following celebrated inequality:
\begin{theorem}[Sun, Theorem 3.1 \cite{Sun}]\label{sun's inequality}
Let $C$ be a smooth projective curve of genus $b\geq 1$ over $k$.
Let $\sE$ be a semi-stable vector bundle of rank $r$ over $C$. One
has the following inequality
$$
\mu_{max}(F_{C}^{*}\sE)-\mu_{min}(F_{C}^{*}\sE)\leq 2(r-1)(b-1).
$$
\end{theorem}
Based on his result, we shall show the following inequality.
\begin{theorem}\label{generalized form of sun's inequality}
Let $C$ be curve as above. The morphism $F_{n}: C\to C, \ n\geq 1$
denotes the $n$-th iterated Frobenius morphism. Let $\sE$ be a
vector bundle of rank $r$ over $C$. We have the following inequality
$$
\mu_{max}(F_{n}^{*}\sE )
 -\mu_{min}(F_{n}^{*}\sE)\leq p^n(\mu_{max}(\sE)-\mu_{min}(\sE))+ 4(r-1)(b-1)\frac{p^n-1}{p-1}.
$$
\end{theorem}
The generalized inequality follows from the following
\begin{proposition}\label{bound of maximal slope}
Let $C$ and $\sE$ be as above. The following inequality holds:
\begin{eqnarray*}
\mu_{max}(F_n^*\sE)\leq p^n\mu_{max}(\sE)+
2(r-1)(b-1)\frac{p^n-1}{p-1}.
\end{eqnarray*}
\end{proposition}
We shall use only Proposition \ref{bound of maximal slope} rather
than Theorem \ref{generalized form of sun's inequality} in our
paper. It is stated for its independent interest as certain
generalization of Sun's inequality. We deduce first Theorem
\ref{generalized form of sun's inequality}
from the above proposition. \\

{\itshape Proof of Theorem \ref{generalized form of sun's
inequality}:} Let
\begin{eqnarray*}
\sE=\mathcal{E}_m\supseteq\mathcal{E}_{m-1}\supseteq\cdots\supseteq
\mathcal{E}_{1}\supseteq\mathcal{E}_{0}=0.
\end{eqnarray*}
be the Harder-Narasimhan filtration of  $\sE$. It is characterized
by two properties, namely the semi-stability of each grading
$\frac{\sE_{i+1}}{\sE_{i}}$ and the strict increase of slopes of
gradings
$$
\mu(\sE_{1})>\mu(\frac{\sE_{2}}{\sE_{1}})>\cdots>\mu(\frac{\sE_{m}}{\sE_{m-1}}).
$$
Let $\sE^\vee$ be the dual vector bundle of $\sE$. It is clear that
the following filtration
$$
\sE^\vee=(\frac{\sE_{m}}{\sE_{0}})^\vee\supset
(\frac{\sE_{m}}{\sE_{1}})^\vee\supset \cdots \supset
(\frac{\sE_{m}}{\sE_{m-1}})^\vee\supset
(\frac{\sE_{m}}{\sE_{m}})^\vee=0
$$
is the Harder-Narasimhan filtration of $\sE^\vee$ by the
characterization. In particular it follows that
\begin{eqnarray*}
 -\mu_{min}(\sE)&=& -\mu(\frac{\sE_{m}}{\sE_{m-1}}) \\
  &= & \mu((\frac{\sE_{m}}{\sE_{m-1}})^\vee) \\
  &=& \mu_{max}(\sE^\vee)  .
\end{eqnarray*}
Thus we apply Proposition \ref{bound of maximal slope} to the vector
bundles $\sE$ and its dual $\sE^\vee$. Hence
\begin{eqnarray*}
\mu_{max}(F_{n}^{*}\sE )
 -\mu_{min}(F_{n}^{*}\sE)&=&\mu_{max}(F_{n}^{*}\sE )
 +\mu_{max}(F_{n}^{*}\sE^\vee)  \\
  &\leq  & p^n\mu_{max}(\sE)+ p^n\mu_{max}(\sE^\vee)+
4(r-1)(b-1)\frac{p^n-1}{p-1} \\
  &=&p^n(\mu_{max}(\sE)-\mu_{min}(\sE))+ 4(r-1)(b-1)\frac{p^n-1}{p-1}.
\end{eqnarray*}
\proofend

So it is left to show Proposition \ref{bound of maximal slope}.\\

{\itshape Proof of Proposition \ref{bound of maximal slope}:} We prove by induction on $n$. \\

{\bf $n=1$ Case.} It is to show
\begin{eqnarray*}\label{n1slop}
\mu_{max}(F_{C}^{*}\sE)\le p\mu_{max}(\sE)+(r-1)(2b-2).
\end{eqnarray*}
We let $\sF\subset F_{C}^{*}\sE$ be the subsheaf with maximal slope.
As above, we let
\begin{eqnarray*}
\sE=\mathcal{E}_m\supseteq\mathcal{E}_{m-1}\supseteq\cdots\supseteq
\mathcal{E}_{1}\supseteq\mathcal{E}_{0}=0.
\end{eqnarray*}
be the Harder-Narasimhan filtration of  $\sE$. It pulls back to a
filtration of $ F_{C}^{*}\sE$:
\begin{eqnarray*}
F_{C}^{*}\sE=F_{C}^{*}\mathcal{E}_m\supseteq
F_{C}^{*}\mathcal{E}_{m-1}\supseteq\cdots\supseteq
F_{C}^{*}\mathcal{E}_{1}\supseteq 0.
\end{eqnarray*}
We consider the natural map $\alpha_1: \mathcal{F}\to
\frac{F_C^{*}\sE_{m}}{F_C^{*}\mathcal{E}_{m-1}}$ which is the
composition of the morphisms
$$
\mathcal{F}\hookrightarrow F_C^{*}\mathcal{E}_m
 \to \frac{F_C^{*}\sE_{m}}{F_C^{*}\mathcal{E}_{m-1}}.
$$
The kernel of $\alpha_1$ is denoted by $\sF_{1}$. So we have a
subsheaf $\frac{\sF}{\sF_{1}}$ of
$\frac{F_C^{*}\sE_{m}}{F_C^{*}\mathcal{E}_{m-1}}$. Since the grading
$\frac{\sE_{m}}{\mathcal{E}_{m-1}}$ is semi-stable, one has after
Theorem \ref{sun's inequality}
$$
\mu_{max}(F_C^{*}(\frac{\sE_{m}}{ \sE_{m-1}}))
 -\mu_{min}( F_C^{*}(\frac{\sE_{m}}{\sE_{m-1}}))\leq 2(r-1)(b-1).
$$
Since furthermore
$$
\mu(\frac{\sF}{\sF_1})\leq \mu_{max}(
\frac{F_C^{*}\sE_{m}}{F_C^{*}\mathcal{E}_{m-1}})= \mu_{max}(
F_C^{*}(\frac{\sE_{m}}{\mathcal{E}_{m-1}}))
$$
and
$$
\mu_{min} (F_C^{*}(\frac{\sE_{m}}{\mathcal{E}_{m-1}}))\leq
\mu(F_C^{*}(\frac{\sE_{m}}{\mathcal{E}_{m-1}}))
=p\mu(\frac{\sE_{m}}{\mathcal{E}_{m-1}})
$$
hold, we get
$$
\mu(\frac{\mathcal{F}}{\sF_1})\leq
2(r-1)(b-1)+p\mu(\frac{\sE_{m}}{\mathcal{E}_{m-1}}) \leq
2(r-1)(b-1)+p\mu_{max}(\sE).
$$
One notes that $\sF_1$ is a coherent subsheaf of
$F_C^{*}\mathcal{E}_{m-1}$. In case that $\sF_1$ is not zero sheaf
one considers further the morphism $\alpha_2: \sF_1\to
\frac{F_C^{*}\sE_{m-1}}{F_C^{*}\mathcal{E}_{m-2}}$ with kernel
$\sF_2$. By keeping on doing this, one obtains a filtration
$$
\sF=\sF_{0}\supset \sF_{1}\supset \cdots \supset \sF_{l}\supset
\sF_{l+1}=0
$$
such that each grading $\frac{\sF_{i}}{\sF_{i+1}}$ is subsheaf of
$F_C^{*}(\frac{\sE_{m-i}}{\mathcal{E}_{m-i-1}})$. One applies then
Theorem \ref{sun's inequality} for each grading
$\frac{\sE_{m-i}}{\mathcal{E}_{m-i-1}}$ and notes that the
inequality
$$
\mu(\frac{\mathcal{E}_{i}}{\mathcal{E}_{i-1}})\leq \mu_{max}(\sE)
$$
holds for all $1\leq i\leq m$. By the same argument as above, one
has
$$
\mu(\frac{\sF_{i}}{\sF_{i+1}})\leq 2(r-1)(b-1)+p\mu_{max}(\sE),\
0\leq i\leq l.
$$
Therefore it follows that
\begin{eqnarray*}
\deg \mathcal{F}&=&\sum_{i=0}^{l} \deg \frac{\sF_i}{\sF_{i+1}}\\
&\leq& \sum
_{i=0}^{l}[2(r-1)(b-1)+p\mu_{max}(\sE)]\rank\frac{\sF_i}{\sF_{i+1}} \\
&=&[2(r-1)(b-1)+p\mu_{max}(\sE)]\rank(\mathcal{F}).
\end{eqnarray*}

{\bf Induction Step.} We show the truth for $n-1$ implies the truth
for $n$. But this is direct. One suffices to notice that the same
argument in the above step applying to the sheaf $F_{n-1}^{*}\sE$
yields the following inequality
$$
\mu_{\max}(F_{n}^{*}\sE)\leq 2(r-1)(b-1)+p\mu_{max}(F_{n-1}^{*}\sE).
$$
By the inductive assumption one has the inequality
$$
\mu_{max}(F_{n-1}^{*}\sE) \leq p^{n-1}\mu_{max} (\sE)+
2(r-1)(b-1)\frac{p^{n-1}-1}{p-1}.
$$
Combining the last two inequalities one gets the claimed inequality
\begin{eqnarray*}
\mu_{max}(F_n^*\sE)\leq p^n\mu_{max}(\sE)+
2(r-1)(b-1)\frac{p^n-1}{p-1}.
\end{eqnarray*}
The proof is completed.

\proofend

{\itshape Proof of (iii):} It follows from Proposition \ref{bound of
maximal slope} for $E^{0,1}$ and (ii).\\

\proofend

In characteristic zero the second Hodge bundle $E^{0,1}$ is
semi-stable and semi-negative. So for any coherent subsheaf
$\sF\subset E^{0,1}$ one has $\mu(\sF)\leq 0$, the same as in the
case (i) of the above theorem. The violation of the semi-stablity in
characteristic $p$ was shown again by the example of Moret-Bailly
\cite{Mo}, in which there is a positive degree sub line bundle in
$E^{0,1}$. It is interesting to ask the stability question about the
Hodge bundles over strata (for example $p$-rank zero stratum or
Ekedahl-Oort strata in general) in the moduli space of algebraic
curves and the moduli space of Abelian varieties over $k$ where the
relative Frobenius degenerates.

\section{Upper Bound of $p$-Rank Zero Locus}
In this section we shall discuss the problem as described in the
introduction, namely the upper bound of $p$-rank zero locus in a
semi-stable family of algebraic curves $f:X\to C$ over $k$. Of
course, in order that our question makes sense we must assume the
generic fiber of $f$ is not of $p$-rank zero. This is the basic
assumption in this section. We denote by $V_{0}(f)$ the proper
subset of $C$ which supports the $p$-rank zero fibers in $f$. The
notation
$|V_{0}(f)|$ means the cardinality of $V_{0}(f)$.\\

We treat first the simplest case, a family of elliptic curves.
\begin{proposition}\label{elliptic curve case}
Let $f: X\to C$ be a semi-stable family of elliptic curves with
nontrivial Kodaira-Spencer map. Then the number of supersingular
elliptic curves in the family $f$ is bounded from above by
$\frac{p-1}{2}(2b-2+s)$ where $b$ is the base curve genus and $s$
the number of singular fibers of $f$.
\end{proposition}
When the base curve is $\P^1$, there are at least four singular
fibers. It follows from the positivity of direct images of relative
differentials over a universal family of Abelian variety in the
dimension one case. So the maximal number of supersingular elliptic
curves is $p-1$, when the base curve is $\P^1$ and the family
degenerates at four points. The six examples of A. Beauville of
semi-stable elliptic curves over $\P^1$ have integral model and
their good reductions at $k$ have exactly $p-1$ supersingular
fibers. It is interesting to ask the converse. Namely, one asks if
the good reduction of Beauville's six examples at a prime $p$ can be
characterized by the minimal number of singular fibers and maximal
number of supersingular fibers.\\

{\itshape Proof:} We consider the morphism $\Phi_{1}:
F_{1}^{*}E^{0,1}\to E^{0,1}$. By Proposition \ref{criteria} (iii),
the supersingular locus $V_{0}(f)$ is the support of the effective
divisor of $E^{0,1}\otimes F_{1}^{*}(E^{0,1})^{-1}$ defined by
$\Phi_{1}$. As shown in the proof of the Arakelov inequality
\ref{Arakelov inequality for families of curves in characteristic
p}, one has inequality
$$
\deg E^{1,0} \leq \frac{1}{2}\deg \Omega_{C}(S).
$$
So it follows that
\begin{eqnarray*}
  |V_{0}(f)| &\leq & \deg(E^{0,1}\otimes F_{1}^{*}(E^{0,1})^{-1}) \\
  &=&  \deg(E^{0,1})-\deg(F_{C}^{*}(E^{0,1})) \\
   &=& (p-1)\deg(E^{1,0})\\
  &\leq &\frac{p-1}{2}(2b-2+s).
\end{eqnarray*}

\proofend

In the above proof the effective divisor of $E^{0,1}\otimes
F_{1}^{*}(E^{0,1})^{-1}=(E^{1,0})^{p-1}$ is classically known as
Hasse locus. When the family is a modular family of elliptic curves
over $k$, it was known that the Hasse locus is reduced. (See the
article \cite{Og} for the multiplicity problem of the Hasse locus of
CY family.) So the first inequality in the above calculation is
indeed an equality for Beauville's six examples. On the other hand,
for these families the second inequality is also an equality because
over $\C$ it is known that the Arakelov inequality reaches equality for them. \\

In the next step we want to bound the cardinality of the
non-ordinary locus, which we denote by $H(f)$, when the family $f$
has an ordinary fiber. One can regard this as the first
generalization of elliptic curve case. Since the $p$-rank zero locus
is contained in the non-ordinary locus, we also obtain an upper
bound for $|V_{0}(f)|$.
\begin{theorem}\label{simple case}
Let $f: X\to C$ be a semi-stable family of algebraic curves over $k$
as above. If there is an ordinary closed fiber in $f$, then we have
the following upper bound for the number of non-ordinary fibers
$|H(f)|$:
$$
|H(f)|\leq 2(p-1)g^2(2b-2+s).
$$
\end{theorem}
{\itshape Proof:} By taking the wedge $g$ power of $\Phi_1$, we
obtain a morphism of invertible sheaves over $C$:
$$
\det(\Phi_1): \bigwedge^{g}F_{1}^{*}(E^{0,1})\to
\bigwedge^{g}E^{0,1}.
$$
It is nontrivial because the family $f$ contains an ordinary fiber
and the morphism $\Phi_1$ is therefore  generically isomorphism by
Proposition \ref{criteria} (i). It must be then injective and we
have an short exact sequence of coherent sheaves over $C$:
$$
0\to \bigwedge^{g}(F_{1}^{*}(E^{0,1}))\stackrel{\det
\Phi_1}{\longrightarrow} \bigwedge^{g}(E^{0,1})\to Q\to 0
$$
where $Q$ is a torsion sheaf. By Proposition \ref{criteria} (i)
again, one has $|H(f)|\leq \deg{Q}$. On the other hand, by the
Arakelov inequality $\deg(E^{1,0})< 2g^{2}(2b-2+s)$ (cf. Theorem
\ref{Arakelov inequality for families of curves in characteristic
p}) it follows that
\begin{eqnarray*}
 |H(f)| &\leq & \deg(Q) \\
  &=&  \deg(E^{0,1})-\deg(F_{1}^{*}(E^{0,1})) \\
   &=& (p-1)\deg(E^{1,0})\\
  &< &2(p-1)g^2(2b-2+s).
\end{eqnarray*}

\proofend

Now we complete our discussions of this problem by considering the
most general case. That is, we are going to provide an upper bound
for $|V_{0}(f)|$ when the family $f$ does not necessarily contain an
ordinary fiber.
\begin{theorem}\label{upper bound of p-rank zero locus}
Let $f: X\to C$ be a semi-stable family of algebraic curves of genus
$g\geq 2$ over $k$ whose Kodaira-Spencer map is nonzero. The
notations is as above. If the generic fiber of $f$ is not of
$p$-rank zero, then the number of $p$-rank zero closed fibers $
|V_{0}(f)|$ in $f$ is bounded strictly from above by the numerical
function
$$
P(p,g,b,s)= [2p^{g}(2g^2-1)+2g(g-1)](2b-2+s).
$$
\end{theorem}
{\itshape Proof:} The proof is in the same line as above. Instead of
considering $\Phi_1$ we need to study $\Phi_g$ in the current case.
Since by assumption the $p$-rank of the generic fiber of $f$ is not
zero, the morphism $\Phi_g$ is nontrivial by Proposition
\ref{criteria} (ii). Thus one obtains the factorization of $\Phi_g$
$$
(F_{g}^{*})E^{0,1} \to \sF \stackrel{\phi}{\rightarrow} \sE \to
E^{0,1}
$$
such that $\phi$ is an isomorphism at the generic point. \\

By taking the wedge product of $\phi$ one has the following short
exact sequence of coherent sheaves over $C$
$$
0\to \det\sF\stackrel{\det\phi}{\longrightarrow} \det \sE \to Q\to 0
$$
where $Q$ is a torsion sheaf. Now one can estimate $|V_{0}(f)|$ as
in the last theorem:
\begin{eqnarray*}
  |V_{0}(f)|&\leq& \deg Q\\
  &= & \deg \sE -\deg \sF \\
  &=&  \deg \sE -\deg(F_{g}^{*}E^{0,1}) +\deg(\ker \Phi) \\
   &= &p^{g}\deg (E^{1,0})+\deg \sE +\deg (\ker \Phi)    \\
   &< & 2p^{g}g^{2}\deg\Omega_{C}(S)+2g(g-1)\deg\Omega_{C}(S)+2p^{g}g(g-1)\deg\Omega_{C}(S)\\
   &&+2(g-1)(b-1)\frac{p^{g}-1}{p-1} \\
&=& [\frac{p^{g+1}-1}{p-1}(g-1)+2p^{g}g(2g-1)+2g(g-1)]\deg\Omega_{C}(S)\\
&&-s(g-1)\frac{p^{g}-1}{p-1}\\
&\leq &
[\frac{p^{g+1}-1}{p-1}(g-1)+2p^{g}g(2g-1)+2g(g-1)](2b-2+s)\\
   &< &[2p^{g}(2g^2-1)+2g(g-1)](2b-2+s).
\end{eqnarray*}
In the first strict inequality above we use the Arakelov inequality
for $E^{1,0}$ together with Theorem \ref{Arakelov inequality for
second Hodge bundle} (ii),(iii). The second strict inequality is
elementary which makes the expression simpler. The whole proof is
completed.

\proofend

\section{Arakelov Inequality and Upper Bound of $p$-Rank Zero Locus for Smooth Families of Abelian Varieties in Characteristic $p$}
In this section we discuss some extensions of previous results for
families of algebraic curves to smooth families of Abelian
varieties. So we let $f:X\to C$ be a smooth family of Abelian
varieties of dimension $g\geq 2$ over $k$. Our basic assumption
about $f$ is as follows:
\begin{assumption}\label{lifting assumption}
Let $f':X'\to C$ be the base change of $f$ under $F_{C}$ (See
Section 3). We assume that $f'$ is $W_2$-liftable. Namely, there is
a smooth family $\tilde f': \tilde X'\to \tilde C$ over $W_{2}(k)$
such that the reduction of $\tilde f'$ at $k$ is $f'$. The specific
choice of lifting of $f'$ is irrelevant to the statements below.
\end{assumption}

We recall first several recent remarkable results due to
Ogus-Vologodsky \cite{OV}.
\begin{proposition}[Ogus-Vologodsky, Theorem 4.14 (3) and Proposition 4.19
\cite{OV}]\label{Higgs semistability} Let $f: X\to C$ be a smooth
family of Abelian varieties over $k$ with assumption \ref{lifting
assumption}. Then the first relative de Rham cohomology
$R^{1}f_{*}^{DR}(\sO_{X})$ of $f$ is a Fontaine module over $C$. Its
grading with respect to the Hodge filtration gives rise to the Higgs
bundle over $C$
$$
(E,\theta)=(E^{1,0}\oplus E^{0,1},\theta^{1,0}\oplus \theta^{0,1}),
$$
which is Higgs semi-stable when $p\geq 4g^2-6g+4$.
\end{proposition}
\begin{theorem}
Let $f: X\to C$ be smooth family of Abelian variety with assumption
\ref{lifting assumption} as above. Assume $p\geq 4g^2-6g+4$. Then
for any coherent subsheaf $\sF$ of $E^{1,0}$ one has inequality of
the slope of $\sF$:
$$
\mu(\sF)\leq \max\{0,b-1\},
$$
where $b$ is genus of $C$. Furthermore, for $\sE$ a coherent
subsheaf $F_{n}^{*}E^{0,1},\ n\geq 0$ one has inequality
$$
\mu(\sE)\leq 2(g-1)(b-1)\frac{p^n-1}{p-1}.
$$
\end{theorem}
{\itshape Proof:} The idea of proof is the same as in characteristic
zero situation (cf. Proposition 1.2 \cite{VZ04}). So let $\sG\otimes
\Omega_C$ be the image of $\sF$ under $\theta^{1,0}$. Then
$\sF\oplus \sG$ forms a Higgs subsheaf of $(E,\theta)$. By the Higgs
semi-stability in Proposition \ref{Higgs semistability}, it follows
that $\deg \sF+\deg \sG\leq \deg E=0$. If $\sG=0$, then one has
$\deg\sF\leq 0$. We assume $\sG\neq 0$ in the following. Hence
\begin{eqnarray*}
 \deg \sF&\leq& \deg \sG+\rank \sG\cdot\deg\Omega_{C}\\
  &\leq & \deg \sG + \rank \sF \cdot\deg\Omega_{C} \\
  &\leq&  -\deg \sF +\rank \sF \cdot\deg\Omega_{C}.
\end{eqnarray*}
It follows that $\mu(\sF)\leq \frac{1}{2}\deg\Omega_{C}=b-1$.\\

Now let $\sE$ be a coherent subsheaf of $E^{0,1}$. One notes that
the second component of Higgs field $\theta^{0,1}$ is simply zero
map. Then $(\sE,0)$ forms a Higgs subsheaf of $(E,\theta)$ and by
Higgs semi-stability \ref{Higgs semistability} $\mu(\sE)\leq 0$. So
the second inequality is proved for $n=0$ case. For generally $\sE$
a coherent subsheaf of $F_{n}^{*}E^{0,1}$, the result follows from
Proposition \ref{bound of maximal slope}. \\

\proofend

The following statement is analogous to that in Theorem \ref{upper
bound of p-rank zero locus}.
\begin{theorem}
Let $f: X\to C$ be a smooth family of algebraic curves of genus
$g\geq 2$ over $k$ with assumption \ref{lifting assumption}. Assume
$p\geq 4g^2-6g+4$. If the generic fiber of $f$ is not of $p$-rank
zero, then the number of $p$-rank zero closed fibers $ |V_{0}(f)|$
in $f$ is bounded from above by the numerical function
$$
Q(p,g,b)= p^g\cdot \max\{0,b-1\}+2(g-1)(b-1)\frac{p^g-1}{p-1}.
$$
\end{theorem}
{\itshape Proof:} The proof is the same as that in Theorem
\ref{upper bound of p-rank zero locus}, except that we replace the
estimates of degrees by those proved in this section. The result can
be easily checked.

\proofend

\end{document}